\documentclass[titlepage,11pt]{article} 
\usepackage{latexsym}
\usepackage{amsfonts}
\usepackage{amsmath}
\newcommand\blackslug{\hbox{\hskip 1pt \vrule width 4pt height 8pt depth 1.5pt
        \hskip 1pt}}
\newcommand\bbox{\hfill \quad \blackslug \bigbreak}

\def\ll{,\ldots,}
\def\cupcup{\cup\cdots\cup}
\newcommand{\vare}{\varepsilon}
\usepackage{tikz}

\title{Erd\H{o}s-Hajnal for graphs with no 5-hole}
\author{Maria Chudnovsky\thanks{Supported by  NSF grant DMS 1763817.}\\
Princeton University, Princeton, NJ 08544
\\
\\
Alex Scott\thanks{Research supported by EPSRC grant EP/V007327/1.}\\
Mathematical Institute, University of Oxford, Oxford OX2 6GG, UK
\\
\\
Paul Seymour\thanks{Supported by AFOSR grant
A9550-19-1-0187, and by NSF grant  DMS-1800053.}\\
Princeton University, Princeton, NJ 08544
\\
\\
Sophie Spirkl\thanks{We acknowledge the support of the Natural Sciences and Engineering Research
Council of Canada (NSERC), [funding reference number RGPIN-2020-03912].
Cette recherche a \'et\'e financ\'ee par le Conseil de recherches en sciences
naturelles et en g\'enie du Canada (CRSNG), [num\'ero de r\'ef\'erence
RGPIN-2020-03912].  }\\
University of Waterloo, Waterloo, Ontario N2L3G1, Canada}

\date{}

\newtheorem{thm}{}[section]

\newcommand{\Proof}{\noindent{\bf Proof.}\ \ }

\begin{document}
\maketitle
\begin{abstract} 
The Erd\H{o}s-Hajnal conjecture says that for every graph $H$ there exists $\tau>0$
such that 
every graph $G$ not containing $H$ as an induced subgraph has a clique or stable set of cardinality at least $|G|^\tau$.
We prove that this is true when $H$ is a cycle of length five. 

We also prove several further results: for instance, that 
if $C$ is a cycle and $H$ is the complement of a forest, 
there exists $\tau>0$ such that every graph $G$ containing neither of $C,H$ as an induced subgraph
has a clique or stable set of cardinality at least $|G|^\tau$.

\end{abstract}

\section{Introduction}

A cornerstone of Ramsey theory is the theorem of Erd\H os and Szekeres \cite{ESz} from the 1930s, that every graph on $n$ vertices
has a clique or stable set of size $\Omega(\log n)$. This order of magnitude cannot be improved, as Erd\H os~\cite{E} showed 
that there are infinitely many graphs $G$ with $\max(\alpha(G),\omega(G))=O(\log(|G|))$, where $\alpha(G)$ and $\omega(G)$ denote 
the cardinalities of (respectively) the largest stable sets and 
cliques in $G$.  Indeed, for most graphs $G$, both $\alpha(G)$ and $\omega(G)$ are of logarithmic size.

The celebrated Erd\H{o}s-Hajnal conjecture asserts that for proper hereditary classes of graphs, the picture is dramatically different.  
We say that a graph $G$ {\em contains} a graph $H$ if some induced subgraph of $G$ is isomorphic to $H$, and $G$ is
{\em $H$-free} otherwise.  Every hereditary class of graphs is defined by its excluded subgraphs (that is, the graphs $F$ such that every graph in the class is $F$-free).
The Erd\H{o}s-Hajnal conjecture~\cite{EH0,EH} asserts that if some graph is excluded, the largest clique or stable set that can be guaranteed jumps from logarithmic to polynomial size:

\begin{thm}\label{EHconj}
{\bf Conjecture: }For every graph $H$, there exists $\tau>0$ such that every $H$-free graph $G$ satisfies
$$\max(\alpha(G),\omega(G))\ge |G|^\tau.$$
\end{thm}

The Erd\H{o}s-Hajnal conjecture is only known for a small family of graphs.  It is trivially true for $H=K_2$; it is true 
for $H=P_4$, the four-vertex path (the $P_4$-free graphs form the well-known class of cographs); and Chudnovsky and 
Safra~\cite{safra} showed that it is true when $H$ is the bull ($P_4$ with an additional vertex adjacent to the two central 
vertices).  It is easy to see that if the conjecture holds for $H$ then it also holds for the complement $\overline H$.  An 
important result of Alon, Pach and Solymosi \cite{aps} shows that if the conjecture holds for $H$ and $H'$ then it also holds for 
the graph obtained by substituting $H'$ into a vertex of $H$.   It follows that the Erd\H{o}s-Hajnal conjecture holds for every 
graph $H$ in the closure of 
$\{K_2, P_4, {\rm bull}\}$ under complements and substitution.  But these are all the graphs (with at least two vertices) for which the conjecture was previously known.

The conjecture holds when $|H|\leq 4$, but is open for three graphs on five vertices: $C_5$, $P_5$ and $\overline{P_5}$.  The five-vertex cycle $C_5$ has been a particularly frustrating open case, and has attracted a good deal of unsuccessful attention (for example, it was highlighted by Erd\H os and Hajnal \cite{EH} and also by Gy\'arf\'as \cite{Gy}).   
So we are happy to report some progress at last: in this paper, we will prove the conjecture for $C_5$, and present a number of other results.

Let us start by noting that the best known bound for a general graph $H$ is due to Erd\H os and Hajnal~\cite{EH}, who showed the following:

\begin{thm}\label{EHpartial}
For every graph $H$, there exists $c>0$ such that 
$$\max(\alpha(G),\omega(G))\ge 2^{c\sqrt{\log |G|}}$$
for every $H$-free graph $G$ with $|G|\geq1$.
\end{thm}

In an earlier paper \cite{c5close}, with Jacob Fox, we improved on \ref{EHpartial} for the graph $H=C_5$:
\begin{thm}\label{mainthm}
There exists $c>0$ such that 
$$\max(\alpha(G),\omega(G))\ge 2^{c\sqrt{\log  |G| \log \log  |G|  }}$$
for every $C_5$-free graph $G$ with $|G|\geq2$.
\end{thm}

Our first result in this paper is the following, proving the Erd\H{o}s-Hajnal conjecture for $C_5$:
\begin{thm}\label{C5thm}
There exists $\tau >0$ such that every $C_5$-free graph $G$ satisfies
$$\max(\alpha(G),\omega(G))\ge |G|^\tau .$$
\end{thm}

The proof of \ref{C5thm} is novel, but the same proof method, with some extra 
twists, yields some other results about the Erd\H{o}s-Hajnal conjecture. It does not seem to show that $P_5$ has the Erd\H os-Hajnal property,
which, with its complement, is the other open case of \ref{EHconj} with $|H|=5$; but it does give other nice things.  In particular, it gives results when certain pairs or small families of induced subgraphs are excluded.

If $\mathcal{H}$ is a set of graphs, $G$ is {\em $\mathcal{H}$-free} if it is $H$-free for each $H\in \mathcal{H}$.
Let $\mathcal{H}$ be a set of graphs (or a single graph); we say that $\mathcal{H}$ {\em has the Erd\H os-Hajnal property}\footnote{Some papers say ``the class of $\mathcal{H}$-free graphs has the Erd\H os-Hajnal property'' in this situation, but here the definition
we give is more convenient.} if there exists
$\tau>0$ such that $\max(\alpha(G),\omega(G))\ge |G|^\tau$ for all $\mathcal{H}$-free graphs.  If $\mathcal{H}=\{H\}$ we simply say that {\em $H$ has the Erd\H os-Hajnal property}.
Thus \ref{EHconj} says that every graph has the Erd\H os-Hajnal property, and \ref{C5thm} says that $C_5$ has the Erd\H os-Hajnal property. Note that if $\mathcal{H}$ has the Erd\H os-Hajnal property then so does the set $\{\overline H:H\in\mathcal H\}$ of complements of members of $\mathcal{H}$.

There has been some recent progress on small sets of graphs with the Erd\H os-Hajnal property.  After partial results by a number of authors (see \cite{lagoutte, hooks, liebenau}), 
the following result was shown in \cite{pure1}:

\begin{thm}\label{treescomp}
If $F$ and $H$ are forests then 
$\{F, \overline{H}\}$ has the Erd\H os-Hajnal property.
\end{thm}

In this paper, we will show that a number of other sets of graphs have
the Erd\H os-Hajnal property.
For instance, let $\widehat{C_5}$ be the graph 
obtained from a cycle $C$ of length five by adding a new vertex with two neighbours in $V(C)$,
adjacent; and in general, let $\overline{H}$ denote the complement graph of a graph $H$.
We will show:
\begin{thm}\label{C5hatthm}
$\{\widehat{C_5},\overline{\widehat{C_5}}\}$ has the Erd\H os-Hajnal property.
\end{thm}
Since $\overline{\widehat{C_5}}$ contains $\overline{P_5}$, this implies that $\{\widehat{C_5}, \overline{P_5}\}$
has the Erd\H os-Hajnal property, strengthening the theorem of \cite{holewithhat} that the set of all ``hole-with-hat'' graphs has the Erd\H os-Hajnal property. 
It also implies the result of Chudnovsky and Safra~\cite{safra} that the bull has the Erd\H os-Hajnal property,
because both $\widehat{C_5},\overline{\widehat{C_5}}$
contain the bull. 

We will show that 
\begin{thm}\label{C6thm}
$\{C_6, \overline{C_6}\}$ has the Erd\H os-Hajnal property.
\end{thm}
and 
\begin{thm}\label{C7thm}
$\{C_7, \overline{C_7}\}$ has the Erd\H os-Hajnal property.
\end{thm}
It would be nice to know if the same is true for $\{C_8, \overline{C_8}\}$, but this remains open.
We will show that one of the forests in \ref{treescomp} can be replaced by a cycle:
\begin{thm}\label{antiforestandhole}
If $C$ is a cycle and $H$ is a forest then $\{C,\overline{H}\}$ has the Erd\H os-Hajnal property.
\end{thm}

We will also show: 
\begin{thm}\label{antiholesandhole}
If $C$ is a cycle and $\ell$ is an integer, the set consisting of $C$ and the complements of all cycles of length at least $\ell$ has the Erd\H os-Hajnal property.
\end{thm}
This strengthens the result of Bonamy, Bousquet and Thomass\'e~\cite{bonamy} that the set consisting of all cycles of length at least $\ell$ and their complements
has the Erd\H os-Hajnal property (see \cite{pure2} for a substantial strengthening of this result).
There are some other more complicated results that we will explain later. 

There are a number of different ways to phrase the Erd\H os-Hajnal conjecture.
Let us define $\kappa(G)=\alpha(G)\omega(G)$.
For a set $\mathcal{H}$ of graphs, the following are equivalent:
\begin{itemize}
\item there exists $\tau>0$ such that every $\mathcal{H}$-free graph $G$ satisfies $\max(\alpha(G),\omega(G))\ge |G|^\tau$;
\item there exists $\tau>0$ such that every $\mathcal{H}$-free graph $G$ contains as an induced subgraph a cograph with at least $|G|^\tau$ vertices (this was implicitly used by Erd\H os and Hajnal \cite{EH});
\item there exists $\tau>0$ such that every $\mathcal{H}$-free graph $G$ contains as an induced subgraph a perfect graph with at least $|G|^\tau$ vertices (this is discussed in \cite{Gy});
\item there exists $\tau>0$ such that every $\mathcal{H}$-free graph $G$ satisfies $\kappa(G)\ge |G|^\tau$.
\end{itemize}
The version using $\kappa$ is sometimes easier to work with, and we will mostly use it below.

An important ingredient in the paper is a lemma about bipartite graphs that we will prove in the next section.  This originates in
a powerful lemma that was proved by Tomon~\cite{tomon}, and developed further by Pach and Tomon~\cite{pach}.
We prove a significant strengthening of Tomon's result, and use it to prove a key lemma that will be used in the proofs of all our main theorems.

The paper is organized as follows. First we prove the strengthening of Tomon's theorem that we need, and then apply it to prove our key lemma; then
we prove \ref{C5thm}; then we extend this approach
to see what else we can obtain, in particular proving the other theorems mentioned above.

Notation throughout is standard. All graphs in this paper are finite and have no loops or parallel edges. We denote by $|G|$ the number of vertices of a graph $G$.
If $X\subseteq V(G)$, $G[X]$ denotes the subgraph of $G$ induced on $X$.
We write $C_k$ for the cycle of length $k$, and $P_k$ for the path with $k$ vertices.
Logarithms are to base two.

\section{A lemma about bipartite graphs}

Let $G$ be a graph.  We say that two sets of vertices $A,B\subseteq V(G)$ are {\em complete} if that are disjoint and 
every element of $A$ is adjacent to every element of $B$, and {\em anticomplete} if they are disjoint and no element of $A$ is 
adjacent to an element of $B$.
We say that a set $\mathcal{H}$ of graphs {\em has the strong Erd\H os-Hajnal property} if there exists $c>0$ such that for every $\mathcal H$-free graph $G$ with at least two vertices there are sets $A,B\subseteq V(G)$ with $|A|,|B|\geq c|G|$ such that the pair $A$, $B$ is either complete or anticomplete.
It is easy to prove that if $\mathcal H$
has the strong Erd\H os-Hajnal property then it has the Erd\H os-Hajnal property (see \cite{semialg,intersection}).  This
approach has been used in a number of papers  to prove the Erd\H os-Hajnal property for various sets $\mathcal H$ (see, for
example,  \cite{bonamy, lagoutte, hooks, pure1,pure2, liebenau}). 

If a finite set $\mathcal H$ of graphs has the strong Erd\H os-Hajnal property then, by considering sparse random graphs, it is easy to see that it is necessary for $\mathcal H$ to contain a forest; and similarly it is necessary for $\mathcal H$ to contain the complement of a forest (see \cite{pure1}).  It follows from \ref{treescomp} that these conditions are also sufficient, and so \ref{treescomp} characterizes finite sets $\mathcal H$ that have the strong Erd\H os-Hajnal property (infinite sets are a different matter: for example the set of all cycles has the strong Erd\H os-Hajnal property, but does not contain a forest).

Tomon \cite{tomon} made the nice observation that there is a similar but weaker property that can also be used to prove the 
Erd\H os-Hajnal property.  Suppose that $\mathcal H$ is a set of graphs and there are $c,k>0$ such that, for every $\mathcal H$-free 
graph $G$ with $|G|\ge 2$, there is some $t=t(G)\geq2$ such that $V(G)$ includes $t$ sets of size at least $c|G|/t^k$ that are pairwise complete 
or pairwise anticomplete (note that the  strong Erd\H os-Hajnal property is the special case where we can always choose
$t=2$).  
We recall that $\kappa(G)=\alpha(G)\omega(G)$; let us write $\kappa(n)$ for the minimum of $\kappa(G)$ over $\mathcal H$-free 
graphs $G$ with $n$ vertices.  It follows that
$\kappa(G)\ge t \kappa(c|G|/t^k)$, and it is easily checked that this implies that $\kappa(n)\ge n^\tau$ for all $n$, provided 
$\tau>0$ is sufficiently small, and so $\mathcal H$ has the strong Erd\H os-Hajnal property.

In order to find the required disjoint sets of vertices, Tomon~\cite{tomon} proved a powerful lemma about bipartite graphs, which was developed further by Pach and Tomon~\cite{pach}.  We will make use of the same idea, but will need to prove a significantly stronger form of the lemma.

Let $G$ be a graph, and let $t,k\ge 0$ where $t$ is an integer. We say $((a_i,B_i)\;:1\le i\le t)$ is a {\em $(t,k)$-comb} in $G$
if:
\begin{itemize}
\item $a_1\ll a_t\in V(G)$ are distinct, and $B_1\ll B_t$ are pairwise disjoint subsets of $V(G)\setminus \{a_1\ll a_t\}$;
\item for $1\le i\le t$, $a_i$ is adjacent to every vertex in $B_i$;
\item for $i,j\in \{1\ll t\}$ with $i\ne j$, $a_i$ has no neighbour in $B_j$; and
\item $B_1\ll B_t$ all have cardinality at least $k$.
\end{itemize}
If $A,B\subseteq V(G)$ are disjoint and $a_1\ll a_t\in A$, and $B_1\ll B_t\subseteq B$, we call this a {\em $(t,k)$-comb in $(A,B)$}.
Our strengthening of Tomon's lemma \cite{tomon} is as follows:

\begin{thm}\label{mainlemma}
Let $G$ be a graph with a bipartition $(A,B)$, such that every vertex in $B$ has a neighbour in $A$; 
and let $\Gamma, \Delta, d>0$ with $d<1$, such that every vertex in $A$ has at most $\Delta$ neighbours in $B$.
Then either:
\begin{itemize}
\item 
for some integer $t\ge 1$, there is a $(t,\Gamma t^{-1/d})$-comb in $(A,B)$; or
\item 
$|B|\le \frac{3^{d+1}}{3/2-(3/2)^{d}}\Gamma^d\Delta^{1-d}$.
\end{itemize}
\end{thm}
\Proof We define a partition of $B$, formed by pairwise disjoint subsets $C_1,C_2,\ldots$ of $B$, defined inductively as follows. 
Let $s\ge 1$, and
suppose that $C_1\ll C_{s-1}$ are defined, and every vertex in $A$ has at most $(2/3)^{s-1}\Delta$ neighbours in $D$, where
$D=B\setminus (C_1\cupcup C_{s-1})$. Choose $a_1,a_2\ll a_k\in A$ with $k$ maximum such that for $1\le i\le k$, there are at
least $(2/3)^{s}\Delta$ vertices in $D$ that are adjacent to $a_i$ and to none of $a_1\ll a_{i-1}$. Let $C_{s}$ be the set of vertices
in $D$ adjacent to one of $a_1\ll a_k$; then from the maximality of $k$, every vertex in $A$ has at most $(2/3)^{s}\Delta$
neighbours in $D\setminus C_{s}$. This completes the inductive definition of $C_1, C_2,\ldots$. Since every vertex in 
$B$
has a neighbour in $A$, it follows that every vertex in $B$ belongs to some $C_s$. 
\\
\\
(1) {\em For all $s\ge 1$, we may assume that $|C_{s}|\le 2^{d+1}(2/3)^{s-sd-1}\Gamma^d\Delta^{1-d}$.}
\\
\\
Let $s\ge 1$, and let $a_1\ll a_k$ be as above (that is, chosen with $k$ maximum such that for $1\le i\le k$, 
there are at
least $(2/3)^{s}\Delta$ vertices in $D=B\setminus (C_1\cupcup C_{s-1})$ that are adjacent to $a_i$ and to none of $a_1\ll a_{i-1}$.) 
For $1\le i\le k$ let $P_i$
be the set of vertices  in $D$ that are adjacent to $a_i$ and to none of $a_1\ll a_{i-1}$; thus each 
$|P_i|\ge (2/3)^{s}\Delta$.
For $1\le i\le k$, let $Q_i$ be the set of vertices in $D\setminus P_i$ adjacent to $a_i$;
thus every vertex in $Q_i$ is adjacent to one of $a_1\ll a_{i-1}$, and 
$$|Q_i| \le (2/3)^{s-1}\Delta - (2/3)^{s}\Delta=(2/3)^{s}\Delta/2$$
since $a_i$ has at most $(2/3)^{s-1}\Delta$ neighbours in $D$.
Inductively, for $i=k, k-1\ll 1$ in turn, we say that $a_i$ is {\em good} if at most
$|P_i|/2$ vertices in $P_i$ are adjacent to a good vertex in $\{a_{i+1}\ll a_k\}$. (Thus $a_k$ is good, if $k>0$.) 
Let $\{a_i\;:i\in I\}$ 
be the set of all good vertices; we claim
that $|I|\ge k/2$. Let $Q$ be the union of the sets $Q_i\;(i\in I)$; then $Q$ has cardinality at most 
$|I|(2/3)^{s}\Delta/2$.
If $i\in \{1\ll k\}\setminus I$, then at least $|P_i|/2\ge (2/3)^{s}\Delta/2$ vertices in $P_i$ belong to $Q$; and so 
$$(k-|I|) (2/3)^{s}\Delta/2\le |Q|\le |I|(2/3)^{s}\Delta/2.$$
Consequently $|I|\ge k/2$.
For each $i\in I$, let $B_i$ be the set of vertices in $P_i$ that are not in $Q$; then $|B_i|\ge  (2/3)^{s}\Delta/2$,
and $((a_i,B_i):i\in I)$ is an $(|I|, (2/3)^{s}\Delta/2)$-comb in $(A,B)$. Let $t=|I|$; so we may assume that either $t=0$, or
$(2/3)^{s}\Delta/2<\Gamma t^{-1/d}$ (since otherwise the theorem holds); and in either case,
$t< (2\Gamma (3/2)^{s}/\Delta)^{d}$. Hence $k\le 2(2\Gamma (3/2)^{s}/\Delta)^{d}$, and 
$$|C_s|\le 2(2\Gamma (3/2)^{s}/\Delta)^{d}\Delta (2/3)^{s-1}.$$
This proves (1).

\bigskip
Now since $d<1$, the sum of $(2/3)^{s-sd-1}$ over all integers $s\ge 1$ equals
$$\frac{(3/2)^{d}}{1-(2/3)^{1-d}},$$
and so 
$$|B|=|C_1|+|C_2|+\cdots \le \frac{2^{d+1}(3/2)^{d}}{1-(2/3)^{1-d}}\Gamma^d \Delta^{1-d}=\frac{3^{d+1}}{3/2-(3/2)^{d}}\Gamma^d\Delta^{1-d}.$$
This proves \ref{mainlemma}.~\bbox

\section{Applying the bipartite lemma}

In this section we use \ref{mainlemma} to prove our key lemma.  Given a vertex $x\in V(G)$, we will apply \ref{mainlemma} to the bipartite graph of edges between $A=N(x)$ and $B=V(G)\setminus(A\cup\{x\})$.  By \ref{mainlemma}, this will either give us a large comb  $((a_i,B_i)\;:1\le i\le t)$, or it will show that $A$ has poor expansion.  In the first case, we try to use the comb either to find $H$ or to find many large sets of vertices that are pairwise complete or pairwise anticomplete; in the second, as there are no edges between $G[A]$ and $G[B\setminus N(A)]$ we can handle them separately.   In both cases, it will be helpful if the set $\{a_1,\dots,a_t\}$ is a stable set: it turns out that we can build this into the key lemma.

Let $\tau>0$. We say that a graph $G$ is {\em $\tau$-critical} if $\kappa(G)<|G|^\tau$, and 
$\kappa(G')\ge |G'|^\tau$ for every induced subgraph $G'$ of $G$ with $G'\ne G$.
The next result is the key lemma that unlocked all the main results in this paper:

\begin{thm}\label{usesparse}
For all $\delta,\vare>0$ with $\vare<1/20$, 
there exists $\tau>0$ with the following property.
Let $G$ be a $\tau$-critical graph, and let $X\subseteq V(G)$
with $|X|\ge \delta|G|$, such that $G[X]$ has maximum degree at most $\vare\delta|G|$.
Then 
there is 
a $(t,\delta|G|/(400\vare t^2))$-comb $((a_i,B_i):1\le i\le t)$ of $G[X]$ such that $t\ge 1/(400\vare)$ and $\{a_1\ll a_t\}$ is stable, and 
there is a vertex $v\in X$ adjacent to $a_1\ll a_t$ and with
no neighbours in $B_1\cupcup B_t$.
\end{thm}
\Proof
Choose $\tau$ with $0<\tau<1$ so small that 
$$\frac{2^{1-1/\tau}}{\delta}+ \left(\vare+\frac{19}{20}\right)\left(\vare\delta\right)^{-\tau}< 1.$$
(This is possible since $\vare<1/20$.)
We claim that $\tau$ satisfies the theorem.

Let $G,X$ be as in the theorem. We may assume that $\kappa(G)\ge 2$.
It follows that $2<|G|^\tau$, and so $|G|>2^{1/\tau}$. 

Let $X_0=X$. Inductively, given a set $X_{i-1}\subseteq X$ with $X_i\ne \emptyset$,
we make the following definitions:
\begin{itemize}
\item Let $v_i\in X_{i-1}$ have maximum degree in $G[X_{i-1}]$.
\item Let $A_i$ be the set of neighbours of $v_i$ in $G[X_{i-1}]$ (possibly $A_i=\emptyset$).
\item Let $C_i\subseteq A_i$ be a stable set with $|C_i|\ge |A_i|^\tau/\omega(G)$. (This exists, since $G$ is $\tau$-critical. Possibly $C_i=\emptyset$, but only if $A_i=\emptyset$.)
\item Let $X_i$ be the set of vertices in $X_{i-1}$ with no neighbour in $\{v_i\}\cup C_i$.
\end{itemize}
\begin{figure}[ht]
\centering

\begin{tikzpicture}[scale=0.7,auto=left]
\draw [rounded corners] (-8, 2) rectangle (2,-2);

\node[inner sep=1.5pt, fill=black,circle,draw] (v1) at ({-7}, {0}) {};

\draw [rounded corners] (-6.2, 1.8) rectangle (-4,-1.8);
\draw [rounded corners] (-5, .5) rectangle (-4.1,-.5);
\draw [rounded corners] (-3.7, 1.8) rectangle (-1,-1.8);
\draw [rounded corners] (-.7, 1.8) rectangle (1.9,-1.8);

\tikzstyle{every node}=[]
\draw (v1) node [left]           {$v_i$};
\node[right] at (-5,0) {$C_i$};
\node[] at (-5.5,0) {$A_i$};
\node[] at (-2.35,0) {$D_i$};
\node[] at (.6,0) {$X_i$};
\node[] at (-3,2.4) {$X_{i-1}$};
\draw[->] (-2.3, 2.4) -- (2, 2.4) {};
\draw[->] (-3.7, 2.4) -- (-8, 2.4) {};
\draw[-] (v1) to (-5.9,0);
\draw[-] (v1) to (-5.9,.3);
\draw[-] (v1) to (-5.9,.6);
\draw[-] (v1) to (-5.9,-.3);
\draw[-] (v1) to (-5.9,-.6);

\draw[-] (-4.25,0) to (-3.2,0);
\draw[-] (-4.25,.2) to (-3.2,.5);
\draw[-] (-4.25,.4) to (-3.2,1);
\draw[-] (-4.25,-.2) to (-3.2,-.5);
\draw[-] (-4.25,-.4) to (-3.2,-1);

\end{tikzpicture}

\caption{Figure for \ref{usesparse}} \label{fig:usesparse}
\end{figure}
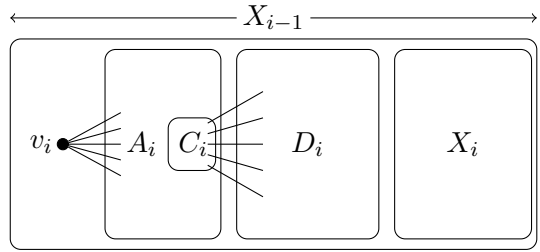
The inductive definition stops when $|X_i|=\emptyset$; let this occur when $i=s$ say. Thus we define a nested sequence of subsets 
$$X=X_0\supseteq X_1\supseteq X_2\supseteq\cdots\supseteq X_s= \emptyset;$$
and also vertices $v_i\in X_{i-1}\setminus X_{i}$ and subsets $A_i,C_i\subseteq X_{i-1}\setminus X_{i}$ for $1\le i\le s$.
Note that there are no edges between $\{v_i\}\cup C_i$ and $X_j$ for $i<j$.

For $1\le i\le s$, let $D_i$ be the set of vertices in $X_{i-1}$ not in $A_i\cup \{v_i\}$,
and with a neighbour in $C_i$. Let $\gamma= \delta/(400\vare)$.
\\
\\
(1) {\em We may assume that $|D_i|\le 19\left(\gamma|A_i|\cdot |X|/\delta\right)^{1/2}$ for $1\le i< s$.}
\\
\\
From the choice of $v_i$, every vertex in $C_i$ has at most $|A_i|$ neighbours in $D_i$.
By \ref{mainlemma} applied to the bipartite graph between $C_i$ and $D_i$, replacing $\Gamma, \Delta, d$ by 
$\gamma|X|/\delta,|A_i|, 1/2$, 
we deduce that either
for some integer $t\ge 1$, there is a $(t,\gamma|X|/(\delta t^2))$-comb in $(C_i,D_i)$,
or
$$|D_i|\le \frac{3^{3/2}}{3/2-(3/2)^{1/2}} (\gamma|X|\cdot|A_i|/\delta)^{1/2}.$$
Suppose the first  holds, and let the comb be $((a_j, B_j):1\le j\le t)$. The sets $B_1\ll B_t$
are pairwise disjoint subsets of $X$, and 
so $t\gamma|X|/(\delta t^2)\le |X|$, that is, $t\ge \gamma/\delta=1/(400\vare)$. Since $|X|\ge \delta|G|$,
it follows that 
$$\gamma|X|/(\delta t^2)\ge \gamma|G|/t^2=\delta |G|/(400\epsilon t^2);$$ 
and therefore in this case
the conclusion of the theorem is true.
So we may assume that the second bullet holds. Since $3^{3/2}/(3/2-(3/2)^{1/2})\le 19$, this proves (1).

\bigskip

For $1\le i\le s$, let $x_i=|A_i|/|X|$.
Since $C_1\cupcup C_{s}$ is stable, and hence has cardinality at most $\alpha(G)$, and 
$$|C_i|\ge \frac{|A_i|^{\tau}}{\omega(G)}= \frac{(x_i|X|)^{\tau}}{\omega(G)}\ge \frac{(x_i\delta|G|)^{\tau}}{\omega(G)}\ge (x_i\delta)^\tau\alpha(G)$$
for each $i$, it follows that 
$\sum_{1\le i\le s} x_i^{\tau}< \delta^{-\tau}$.

Now $X$ is partitioned into the sets $\{v_i\}\;(1\le i\le s)$, $A_i\;(1\le i\le s)$ and $D_i\;(1\le i\le s)$, and so
$$\sum_{1\le i\le s}(1+|A_i|+|D_i|)=|X|,$$
that is,
$$\frac{s}{|X|}+\sum_{1\le i\le s}\frac{|A_i|}{|X|} + \sum_{1\le i\le s}\frac{|D_i|}{|X|}=1.$$
We will bound these three terms separately.

First, since $\{v_1\ll v_s\}$ is stable, it follows that 
$$\frac{s}{|X|}\le \frac{\alpha(G)}{|X|}\le \frac{|G|^\tau}{|X|}\le \frac{|G|^{\tau-1}}{\delta},$$
and since $|G|^{\tau-1}\le 2^{1-1/\tau}$ (because $|G|\ge 2^{1/\tau}$),
it follows that $s/|X|< 2^{1-1/\tau}/\delta$.

Second, 
$$\sum_{1\le i\le s}\frac{|A_i|}{|X|}=\sum_{1\le i\le s}x_i=\sum_{1\le i\le s}x_i^{\tau}x_i^{1-\tau}\le \sum_{1\le i\le s}x_i^{\tau}\vare^{1-\tau}\le \vare(\vare\delta)^{-\tau}$$
since $x_i=|A_i|/|X|\le \vare\delta|G|/|X|\le \vare$.

Third, 
$$\sum_{1\le i\le s}\frac{|D_i|}{|X|}\le  19 \left(\frac{\gamma}{\delta}\right)^{1/2}\sum_{1\le i\le s}x_i^{1/2}= \frac{19}{20}\vare^{-1/2}\sum_{1\le i\le s}x_i^{1/2}$$ 
by (1) and the definition of $\gamma$; 
and 
$$\sum_{1\le i\le s}x_i^{1/2}=\sum_{1\le i\le s}x_i^{\tau}x_i^{1/2-\tau}\le \sum_{1\le i\le s}x_i^{\tau}\vare^{1/2-\tau}\le 
\vare^{1/2}(\vare\delta)^{-\tau}.$$
Consequently 
$$\sum_{1\le i\le s}\frac{|D_i|}{|X|} \le \frac{19}{20}(\vare\delta)^{-\tau}.$$
Summing, we deduce that 
$$\frac{2^{1-1/\tau}}{\delta}+ \left(\vare+\frac{19}{20}\right)(\vare\delta)^{-\tau}\ge 1,$$
contrary to the choice of $\tau$. This proves \ref{usesparse}.~\bbox

This gives us a $(t,\delta|G|/(400\vare t^2))$-comb. The $t^2$ in the denominator comes from applying \ref{mainlemma} with 
$d=1/2$; as was observed by Pach and Tomon~\cite{pach}, we could apply \ref{mainlemma} with $d=1/k$, for any real $k>1$, and produced a comb with $t^{k}$ in the denominator, 
but there is no gain for us in the applications.

\section{The simplest application: excluding $C_5$}

In this section we prove \ref{C5thm}. This is implied by each of several stronger results later in the paper, 
but since the $C_5$ result is of great interest,
and the argument for $C_5$ is easier than the material to come later (which will require additional ideas), we give a separate proof.
We will need
R\"odl's theorem~\cite{rodl}:
\begin{thm}\label{prerodl}
For every graph $H$ and all $\vare>0$, there exists $\delta>0$ such that for every $H$-free graph $G$, there exists $X\subseteq V(G)$
with $|X|\ge \delta|G|$, such that one of $G[X], \overline{G}[X]$ has  at most $\vare |X|(|X|-1)$ edges.
\end{thm}

We also need the following:
\begin{thm}\label{density}
Let $G$ be a graph with at most $\vare|G|(|G|-1)/2$ edges; then for every integer $m\ge 0$ with $m\le (|G|+1)/2$, there exists
$X\subseteq V(G)$ with $|X|=m$ such that $G[X]$ has maximum degree less than $\vare (m-1)$.
\end{thm}
\Proof
By averaging over all subsets $Y$ of $V(G)$ with cardinality $2m-1$, it follows that there exists such a set $Y$ where $G[Y]$
has at most $\vare(2m-1)(m-1)<2\vare m(m-1)$ edges. Thus fewer than $m$ vertices in $Y$ have at least $\vare (m-1)$ neighbours in $Y$.
Hence there exists $X\subseteq Y$ with $|X|=m$ such that $G[X]$ has maximum degree less than $\vare (m-1)$. This proves \ref{density}.~\bbox

We deduce a slight but convenient strengthening of \ref{prerodl}, the following (this is well-known, but we include the proof for completeness):
\begin{thm}\label{rodl}
For every graph $H$ and all $\vare>0$, there exists $\delta>0$ such that for every $H$-free graph $G$, there exists $X\subseteq V(G)$
with $|X|\ge \delta|G|$, such that one of $G[X], \overline{G}[X]$ has maximum degree at most $\vare\delta|G|$.
\end{thm}
\Proof
Let $\vare'=\vare/2$. By \ref{prerodl} there exists $\delta'>0$ such that for every $H$-free graph $G$, there exists $Z\subseteq V(G)$
with $|Z|\ge \delta'|G|$, such that one of $G[Z], \overline{G}[Z]$ has  at most $\vare' |Z|(|Z|-1)$ edges. Let $\delta=\delta'/2$; we claim
that $\delta$ satisfies the theorem. Thus, let $G$ be $H$-free. By the choice of $\delta'$, there exists $Z\subseteq V(G)$
with $|Z|\ge \delta'|G|$, such that one of $G[Z], \overline{G}[Z]$ has  at most $\vare' |Z|(|Z|-1)=\vare |Z|(|Z|-1)/2$ edges. By replacing $G$
by its complement if necessary, we may assume the first. Let $m=\lceil \delta|G|\rceil$; then 
$$|Z|\ge \lceil\delta'|G|\rceil\ge 2m-1.$$
By \ref{density} applied to $G[Z]$, there exists
$X\subseteq Z$ with $|X|=m$ such that $G[X]$ has maximum degree less than $\vare (m-1)\le \vare\delta|G|$. This proves \ref{rodl}.~\bbox

If $X\subseteq V(G)$, we sometimes write $\alpha(X)$ for $\alpha(G[X])$ and so on.
Now we can prove the main result of this section, which we restate:
\begin{thm}\label{C5thm2}
$C_5$ has the Erd\H os-Hajnal property.
\end{thm}
\Proof
Choose $\vare$ with $0<\vare<1/400$, and choose
$\delta$ satisfying \ref{rodl} with $H=C_5$.
Let $\tau>0$ satisfy \ref{usesparse}. Every positive number smaller than $\tau$ also satisfies \ref{usesparse}, and 
since $400\vare<1$, by reducing $\tau$ we may assume that $(400\vare)^{2-1/\tau}>400\vare/\delta$.
We will show that $\kappa(G)\geq|G|^\tau$ for every $C_5$-free graph $G$.  By the remarks in the introduction, this is equivalent to showing that $C_5$ has the Erd\H os-Hajnal property.

Suppose that there is a  $C_5$-free graph $G$ with $\kappa(G)<|G|^\tau$, and choose $G$ minimal;
then $G$ is $\tau$-critical.
By \ref{rodl} there exists $X\subseteq V(G)$
with $|X|\ge \delta|G|$, such that one of $G[X], \overline{G}[X]$ has maximum degree at most $\vare\delta|G|$.
By replacing $G$ with its complement if necessary (this is legitimate since $\overline{G}$ is also $C_5$-free and $\tau$-critical)
we may assume that $G[X]$ has maximum degree at most $\vare\delta|G|$. By \ref{usesparse} and the choice of $\tau$, 
there is
a $(t,\delta|G|/(400\vare t^2))$-comb $((a_i,B_i):1\le i\le t)$ of $G[X]$ such that $t\ge 1/(400\vare)$ and $\{a_1\ll a_t\}$ is stable, and
there is a vertex $v\in X$ adjacent to $a_1\ll a_t$ and with
no neighbours in $B_1\cupcup B_t$.

If there exist $i,j$ with $1\le i<j\le t$ such that some vertex $b_i\in B_i$ has a neighbour $b_j\in B_j$, then the subgraph
induced on $\{b_1,b_2,a_1,a_2,v\}$ is isomorphic to $C_5$, a contradiction. So the sets $B_1\ll B_t$ are pairwise anticomplete. 
Since $G$ is $\tau$-critical, it follows that $\kappa(B_i)\ge |B_i|^\tau$ for each $i$, and since $\kappa(B_i)\le \alpha(B_i)\omega(G)$,
we have 
$$\alpha(B_i)\ge |B_i|^\tau/\omega(G)\ge (\delta|G|/(400\vare t^2))^\tau/\omega(G).$$ 
Since  $B_1\ll B_t$ are pairwise anticomplete, it follows that
$$\alpha(G)\ge \sum_{1\le i\le t}\alpha(B_i)\ge t (\delta|G|/(400\vare t^2))^\tau/\omega(G),$$
and so $\kappa(G)\ge t(\delta|G|/(400\vare t^2))^\tau$. Since $\kappa(G)<|G|^\tau$, it follows that
$400\vare/\delta \ge t^{1/\tau-2}$. But $t\ge 1/(400\vare)$, and $\tau<1/2$, and so
$400\vare/\delta \ge (400\vare)^{2-1/\tau}$,
contrary to the choice of $\tau$. This proves~\ref{C5thm2}.~\bbox

\section{Blockades}

Next we will add some refinements to the proof of \ref{C5thm}, but first
let us set up some more terminology. Let $G$ be a graph.
A {\em pure pair} in $G$ is a pair $A,B$ of disjoint subsets of $V(G)$ such that $A$ is either complete or anticomplete to $B$.
A {\em blockade} $\mathcal{B}$ in $G$ is a sequence $(B_1\ll B_t)$ of pairwise disjoint subsets
of $V(G)$ called {\em blocks}. (In
this paper the order of the blocks $B_1\ll B_t$ in the sequence will not matter.) We denote $B_1\cupcup B_t$ by
$V(\mathcal{B})$. The {\em length} of a blockade is the number of blocks, and its {\em width} is the minimum cardinality of a block.

A blockade $\mathcal{B}=(B_1\ll B_t)$ in $G$ is {\em pure} if $(B_i, B_j)$
is a pure pair for all $i,j$ with $1\le i<j\le t$.
Let $P$ be the graph with vertex
set $\{1\ll t\}$, in which $i,j$ are adjacent if $B_i$ is complete to $B_j$. We say $P$ is the {\em pattern} of the
pure blockade $\mathcal{B}$. A {\em cograph} is a $P_4$-free graph. Every cograph $P$ with more than one vertex admits a pure pair
$(A,B)$ with $A,B\ne \emptyset$ and with $A\cup B=V(P)$.

We need:
\begin{thm}\label{cographpattern}
Let $\mathcal{B}=(B_1\ll B_t)$ be a pure blockade with a cograph pattern.
Then
$$\kappa(B_1\cupcup B_t)\ge \sum_{1\le i\le t}\kappa(B_i).$$
\end{thm}
\Proof
We proceed by induction on $t$. If $t=1$ the claim is true, so we assume $t>1$. Hence
there is a partition $(I,J)$ of
$\{1\ll t\}$, with $I,J\ne \emptyset$, such that either $B_i$ is complete to $B_j$ for all $i\in I$ and $j\in J$, or
$B_i$ is anticomplete to $B_j$ for all $i\in I$ and $j\in J$. We may assume the second by replacing $G$ by its complement if necessary.
Let $V=V(\mathcal{B})$, and
$U=\bigcup_{i\in I}B_i$, and $W=\bigcup_{j\in J} B_j$. Thus $U$ is anticomplete to $W$.
From the inductive hypothesis, $\kappa(U)\ge   \sum_{i\in I}\kappa(B_i)$ and $\kappa(W)\ge  \sum_{j\in J}\kappa(B_j)$. But
$$\kappa(V)=\alpha(V)\omega(V)= (\alpha(U)+\alpha(W))\omega(V)\ge \alpha(U)\omega(U)+\alpha(W)\omega(W)=\kappa(U)+\kappa(W),$$
and the result follows. This proves \ref{cographpattern}.~\bbox

The following is a slight extension of an idea of Pach and Tomon~\cite{pach} (which they called the ``quasi-Erd\H{o}s-Hajnal property''):
\begin{thm}\label{parties}
Let $\tau>0$, and suppose that $G$ is $\tau$-critical.
Then for every integer $t>0$, there is no pure blockade in $G$ with a cograph pattern,
of length $t$ and width at least $|G|t^{-1/\tau}$, such that $B_i\ne V(G)$ for each $i$.
\end{thm}
\Proof
Suppose that $\mathcal{B}=(B_1\ll B_t)$ is such a blockade. Since $G$ is $\tau$-critical, $\kappa(B_i)\ge |B_i|^\tau\ge |G|^\tau/t$ for each $i$,
and so by \ref{cographpattern},
$$\kappa(G)\ge \kappa(B_1\cupcup B_t)\ge \sum_{1\le i\le t}\kappa(B_i)\ge |G|^\tau,$$
a contradiction. This proves \ref{parties}.~\bbox

\section{Forests and their complements}

The proof of \ref{C5thm2} can be developed to give more. 
We have two ways to do so, and in this section we explain the first.

If $H$ is a graph, we wish to augment it in two ways. Let the vertices of $H$ be $\{b_1\ll b_k\}$, and add $k+1$ 
new vertices $a_1\ll a_k,v$ to $V(H)$, where $a_i$ is adjacent to $b_i$ for $1\le i\le k$, and $v$ is adjacent to
$a_1\ll a_k$, and there are no other edges. Let the graph we obtain be $H'$. 
We call $H'$ a {\em star-expansion} of $H$.

In this section we will prove:
\begin{thm}\label{treeext}
Let $H$ be a forest. Let $H_1$ be the star-expansion of $H$, and let $H_2$
be the star-expansion of $\overline{H}$. Then 
$$\{H_1,H_2,\overline{H_1},\overline{H_2}\}$$
has the Erd\H os-Hajnal property.
\end{thm}

\begin{figure}[ht]
\centering

\begin{tikzpicture}[scale=0.7,auto=left]
\tikzstyle{every node}=[inner sep=1.5pt, fill=black,circle,draw]

\node (v) at (-4,5) {};
\node (a1) at (-7,3) {};
\node (a2) at (-5,3) {};
\node (a3) at (-3,3) {};
\node (a4) at (-1,3) {};
\node (b1) at (-7,1) {};
\node (b2) at (-5,1) {};
\node (b3) at (-3,1) {};
\node (b4) at (-1,1) {};

\foreach \from/\to in {v/a1,v/a2,v/a3,v/a4,a1/b1,a2/b2,a3/b3,a4/b4,b1/b2,b2/b3, b3/b4}
\draw [-] (\from) -- (\to);

\end{tikzpicture}

\caption{The star-expansion of $P_4$.} \label{fig:P4expansion}
\end{figure}
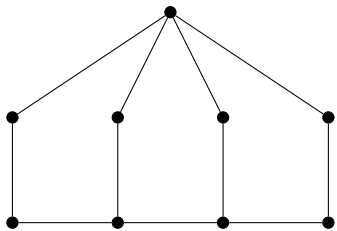

\ref{treeext} is particularly nice when $H=P_4$, since $P_4$ is isomorphic to its complement, and so we only have to exclude two 
graphs instead of four. We obtain:
\begin{thm}\label{twoleaves}
Let $H$ be the graph of figure~\ref{fig:P4expansion}; then $\{H,\overline{H}\}$ has the Erd\H os-Hajnal property.
\end{thm}

This contains \ref{C5thm}, 
\ref{C6thm} and \ref{C7thm}, because the graph of figure \ref{fig:P4expansion}
contain $C_5$, $C_6$ and $C_7$. (The approach via \ref{treeext} does not work for $C_8,\overline{C_8}$,
because there is no forest $H$ such that the star-expansion
of $\overline{H}$ contains one of $C_8,\overline{C_8}$.)

If $\mathcal{B}=(B_1\ll B_t)$ is a blockade in $G$, we say an induced subgraph $H$ of $G$ is {\em $\mathcal{B}$-rainbow} if $V(H)\subseteq V(\mathcal{B})$ and
$|B_i \cap V(H)|\le 1$ for $1\le i\le t$.
To prove \ref{treeext} we need the following theorem of \cite{pure1}:

\begin{thm}\label{pure1}
For every forest $H$, there exist $d>0$ and an integer $K$ with the following property. Let $G$ be a graph 
with a blockade $\mathcal{B}$ of length at least $K$, and let $W$ be the width of $\mathcal{B}$. If 
every vertex of $G$ has degree less than $W/d$, and
there is no anticomplete pair $A,B\subseteq V(G)$ with $|A|,|B|\ge W/d$,
then
there is a $\mathcal{B}$-rainbow copy of $H$ in $G$.
\end{thm}

We used this in \cite{pure1} to deduce that for every forest $H$, the set $\{H,\overline{H}\}$ has the Erd\H os-Hajnal property.
We see that \ref{treeext} (applied to the forest $H$) will be an
extension of that result, since the four graphs of \ref{treeext} all contain one of $H, \overline{H}$.

To prove \ref{treeext}, we need to bootstrap \ref{pure1} into something stronger, and we do so in several stages. 
We will use a strengthening of \ref{prerodl}, 
due to Nikiforov~\cite{nikiforov}, the following:

\begin{thm}\label{prefoxrodl}
For all $\vare>0$ and every graph $H$ on $h$ vertices, there exist $\gamma, \delta>0$ such that
if $G$ is a graph containing
fewer than $\gamma|G|^h$ induced labelled copies of $H$, then there exists $X\subseteq V(G)$ with $|X|\ge \delta|G|$
such that one of $G[X], \overline{G}[X]$ has at most $\vare|X|(|X|-1)$ edges.
\end{thm}

Applying \ref{density} as before yields:

\begin{thm}\label{foxrodl}
For all $\vare>0$ and every graph $H$ on $h$ vertices, there exist $\gamma, \delta>0$ such that
if $G$ is a graph containing
fewer than $\gamma|G|^h$ induced labelled copies of $H$, then there exists $X\subseteq V(G)$ with $|X|\ge \delta|G|$
such that one of $G[X], \overline{G}[X]$ has maximum degree at most $\vare\delta|G|$.
\end{thm}

Let us prove a version of \ref{pure1} without the sparsity hypothesis:

\begin{thm}\label{rainbowforest}
For every forest $H$, there exist $d>0$ and an integer $K$, such that,
for every graph $G$ with a blockade $\mathcal{B}$ of length at least $K$, if 
there is no pure pair $A,B\subseteq V(G)$ with $|A|,|B|\ge W/d$, where $W$ is the width
of $\mathcal{B}$,
then
there is a $\mathcal{B}$-rainbow copy of one of $H,\overline{H}$ in $G$.
\end{thm}
\Proof
Choose $d', K'$ to satisfy \ref{pure1} (with $d,K$ replaced by $d', K'$). Let 
$\vare\le 1/(2d'K')$ with $\vare>0$,
and choose $\gamma, \delta>0$ to satisfy \ref{foxrodl}. Choose $K\ge 2K'/\delta$, and such that
$(1-h/K)^h>1-\gamma$.
Choose $d$ such that
$d\ge d'K'/(\delta K-K')$. We claim that $K,d$ satisfy the theorem. 
\\
\\
(1) {$\delta K/K'-1\ge \max\left(\vare \delta d'K,d'/d\right)$.}
\\
\\
To see that $\delta K/K'-1 \ge \vare \delta d'K$, observe that $\delta K/(2K')\ge 1$, and 
$\delta K/(2K')\ge \vare \delta d'K$. The second part, that $\delta K/K'-1 \ge d'/d$, is true from the choice of $d$.
This proves (1).

\bigskip

Let $G$ be a graph with a blockade $\mathcal{B}=(B_1\ll B_K)$ and of width $W$. We may assume that
$|B_i|=W$ for each $i$, and so $|V|=KW$, where $V=B_1\cupcup B_K$.
We assume that there is no $\mathcal{B}$-rainbow copy of $H$. But the number of sequences $(v_1\ll v_h)$
with $v_1\ll v_h\in V$, such that $v_1\ll v_h$ all belong to different blocks of the blockade, is 
$$|V|(|V|-W)(|V|-2W)\cdots (|V|-(h-1)W)> (1-h/K)^h|V|^h\ge (1-\gamma)|V|^h,$$
and since none of them induce a $\mathcal{B}$-rainbow copy of $H$, it follows that the number of induced labelled copies of $H$
in $G[V]$ is less than $|V|^h-(1-\gamma)|V|^h=\gamma|V|^h$.
By \ref{foxrodl} applied to $G[V]$, there exists $X\subseteq V$ with $|X|\ge \delta KW$, such that
one of $G[X], \overline{G}[X]$ has maximum degree less than $\vare\delta KW$; and by replacing $G$ by $\overline{G}$ if
necessary, we may assume that $G[X]$ has maximum degree less than $\vare\delta KW$.
By (1), there exists a real number $W'$ such that
$$\frac{\delta K}{K'}-1\ge \frac{W'}{W}\ge \max\left(\vare \delta d'K,\frac{d'}{d}\right).$$

The sets $B_1\cap X\ll B_K\cap X$ each have cardinality at most $W$, but their union has cardinality at least $\delta KW$.
Let us choose pairwise disjoint subsets $I_1\ll I_t$ of $\{1\ll K\}$, with $t$ maximum such that $|B_h'|\ge W'$
for $1\le h\le t$, where $B_h'=\bigcup_{i\in I_h}B_i\cap X$.
We may assume that $I_1\ll I_{t}$ are minimal with this property, and so $|B_h'|\le W'+W$
for $1\le h\le t$. From the maximality of $t$, 
$$\sum(|B_i\cap X|:i\in \{1\ll t\}\setminus I_1\cupcup I_t)<W';$$
and so $\delta KW\le |X|\le t(W'+W)+W'$. Since $\delta K/K'-1\ge W'/W$ it follows that $t\ge K'$.

Let $\mathcal{B}'$ be the blockade $(B_1'\ll B_{K'}')$; it has width at least $W'$. 
\\
\\
(2) {\em There is a $\mathcal{B}'$-rainbow
copy of $H$.}
\\
\\
Suppose not. Let 
$G'=G[B_1'\cupcup B_{K'}']$. By \ref{pure1} applied to $G'$, 
either
\begin{itemize}
\item some vertex of $G'$ has degree at least $W'/d'$; or
\item there is an anticomplete pair $A,B\subseteq V(G')$ with $|A|,|B|\ge W'/d'$.
\end{itemize}
But the first does not hold, since $G'$ has maximum degree less than $\vare\delta KW\le W'/d'$; and the second
does not hold, since $W'/d'\ge W/d$ and there is no pure pair
$(A,B)$ in $G$ with $|A|,|B|\ge  W/d$. This proves (2).

\bigskip
The $\mathcal{B}'$-rainbow copy of $H$ in (2) is also $\mathcal{B}$-rainbow. This proves \ref{rainbowforest}.~\bbox

\begin{thm}\label{blockparty}
For every forest $H$, there exist an integer $d>0$, such that,
for every integer $s\ge 1$ and every graph $G$, the following holds. Let $D=2^{s-1}d^{2s-1}$, and let $\mathcal{B}$ be a 
blockade in $G$ of length $D$. Then 
either
\begin{itemize}
\item
$G$ admits a pure blockade $\mathcal{A}$ with a cograph pattern, of length $2^s$ and width at least $W/D$, 
where $W$ is the width of $\mathcal{B}$; or
\item there is a $\mathcal{B}$-rainbow copy of one of $H,\overline{H}$ in $G$.
\end{itemize}
\end{thm}
\Proof
Choose $K,d$ to satisfy \ref{rainbowforest}. Then any pair of numbers $K',d'$ with $K'\ge K$ and $d'\ge d$ also 
satisfy \ref{rainbowforest}, so by increasing $K$ or $d$ if necessary, we may assume that $K=d$. We claim that $d$ satisfies \ref{blockparty}.
This is true if $s=1$, from the choice of $d$, and so we assume it is true for some $s\ge 1$ and prove it for $s+1$.

Let $D=2^{s}d^{2s+1}$, and let $G$ be a graph with a blockade $\mathcal{B}=(B_1\ll B_D)$ of width $W$. Partition $\{1\ll D\}$
into $d$ sets of cardinality $D/d$, say $I_1\ll I_d$. Let $B_h'=\bigcup_{i\in I_h}B_i$ for $1\le i\le d$; then 
$\mathcal{B}'=(B_1'\ll B_d')$ is a blockade of length $d$ and width $WD/d$. Let $G'=G[B_1\cupcup B_D]$. 
We may assume there is no $\mathcal{B}'$-rainbow copy of 
$H$ or of $\overline{H}$ in $G'$; so from the choice of $d$, there is a pure pair $(A,B)$ of $G'$ with $|A|,|B|\ge WD/d^2$.

Let $W'=W/(2d^2)$, and $D'=2^{s-1}d^{2s-1}$. 
Let $p$ be the number of $i\in \{1\ll D\}$ such that $|A\cap B_i|\ge W'$. Then 
$pW+DW'\ge |A|\ge WD/d^2$, and so $p\ge D/(2d^2)=D'$. Let $\mathcal{C}$ be the blockade formed by the $D'$
largest sets of the form $A\cap B_i$; then $\mathcal{C}$ has width at least $W'$, and we may assume that
there is no $\mathcal{C}$-rainbow copy of 
$H$ or of $\overline{H}$. Thus the inductive hypothesis, applied to the blockade $\mathcal{C}$ of $G[A]$
implies that $G[A]$ admits a 
pure blockade with a cograph pattern, of width at least $W'/D'=W/D$ and length $2^s$; and similarly so does $G[B]$.
But then combining these gives a pure blockade in $G$ with a cograph pattern, of width at least $W/D$ and length $2^{s+1}$.
This proves \ref{blockparty}.~\bbox

Now finally we can prove \ref{treeext}, which we restate:
\begin{thm}\label{treeext2}
Let $H$ be a forest. Let $H_1$ be the star-expansion of $H$, let $H_2$
be the star-expansion of $\overline{H}$, and let 
$$\mathcal{H}=\{H_1,H_2,\overline{H_1},\overline{H_2}\}.$$
Then $\mathcal{H}$ has the Erd\H os-Hajnal property.
\end{thm}
\Proof
Much of the proof is the same as for \ref{C5thm2}.
Let $d$ satisfy \ref{blockparty}.
Choose $\vare>0$ with $\vare <1/(400d)$,
choose $\delta$ to satisfy \ref{rodl} with $H=H_1$, and let $\gamma=\delta/(400\vare)$.
Choose $\tau>0$ satisfying \ref{usesparse}, such that $1/\tau>3+6\log_2(d)$, and such that
$2^{q}d^{2q+1}<1/(400\vare)$  where
$$q=\frac{\log_2(d)-\log_2(\gamma)}{1/\tau-3-6\log_2(d)}.$$
(We can satisfy the last condition since by making $\tau$ sufficiently small we can make $q$ arbitrarily close to $0$,
and hence make $2^{q}d^{2q+1}$ arbitrarily close to $d< 1/(400\vare)$.)

As in the proof of \ref{C5thm2}, we may assume that there is a $\tau$-critical $\mathcal{H}$-free graph $G$,
and there exists $X\subseteq V(G)$
with $|X|\ge \delta|G|$, such that $G[X]$ has maximum degree at most $\vare\delta|G|$.

By \ref{usesparse} and the choice of $\tau$,
there is
a $(t,\gamma|G|/t^2)$-comb $((a_i,B_i):1\le i\le t)$ of $G[X]$ such that $t\ge 1/(400\vare)$ and $\{a_1\ll a_t\}$ is stable, and
there is a vertex $v\in X$ adjacent to $a_1\ll a_t$ and with
no neighbours in $B_1\cupcup B_t$.
Let $\mathcal{B}=(B_1\ll B_t)$. 
\\
\\
(1) {\em There is a $\mathcal{B}$-rainbow copy of $H$ or of $\overline{H}$.}
\\
\\
Suppose not. Choose an integer $s$ maximum such that $D_s\le t$, where $D_s=2^{s-1}d^{2s-1}$. 
Thus $D_{s+1}>t$. Since $D_{q+1}\le 1/(400\vare)\le t$, it follows that $s>q$.

By \ref{blockparty}, 
$G$ admits a pure blockade $\mathcal{A}$ with a cograph pattern, of width at least $\gamma|G|/(t^2D_s)$ and length $2^s$.
By \ref{parties}, 
$\gamma|G|/(t^2D_s) < |G|(2^{s})^{-1/\tau}$, that is,
$\gamma<t^2D_s2^{-s/\tau}$. The maximality of $s$ implies that $2^{s}d^{2s+1}\ge t$, and so,
substituting for $t$ and for $D_s$, 
we obtain
$$\gamma<2^{2s}d^{4s+2}2^{s-1}d^{2s-1}2^{-s/\tau}.$$
It follows that $\log_2(\gamma) + s/\tau-3s+1<(6s+1)\log_2(d)$,
and so
$$\left(\frac{1}{\tau}-3-6\log_2(d)\right)s<\log_2(d)-\log_2(\gamma).$$
Hence $s<q$, a contradiction.
This proves (1).

\bigskip

But now the result follows as in \ref{C5thm2}.
This proves \ref{treeext}.~\bbox

\section{Excluding a forest complement}

If $H$ is a forest, then since two of the four graphs of \ref{treeext} contain $\overline{H}$, it follows that the set
consisting of $\overline{H}$ and
the remaining two graphs in \ref{treeext} has the Erd\H os-Hajnal property. But we can do better than this: it is sufficient just to exclude
one of the remaining two, as we show in this section. This is proved by
a slight variation in the proof of \ref{treeext}.

We will need the following theorem of~\cite{pure1} (it is a consequence of \ref{pure1}):
\begin{thm}\label{pure1weak}
For every forest $H$, there exists $\vare>0$ such that if a graph $G$ with $|G|>1$ has maximum degree less than $\vare|G|$, and has
no anticomplete pair of sets $A,B\subseteq V(G)$ with $|A|,|B|\ge \vare|G|$, then $G$ contains $H$.
\end{thm}

We use this to prove:
\begin{thm}\label{antiforest}
Let $H$ be a forest, and let $H'$ be the star-expansion of $H$. Then $\mathcal{H}=\{\overline{H}, H'\}$
has the Erd\H os-Hajnal property.
\end{thm}
\Proof
We define $d,\vare,\delta,\tau$ and the rest, exactly as in the proof of \ref{treeext}, except we choose $\vare$ satisfying
\ref{pure1weak} as well as the other conditions, and choose $\tau$ such that  $\vare\delta>2^{-1/\tau}$ as well as the 
other conditions.

As before, we may assume that there is a $\tau$-critical $\mathcal{H}$-free graph $G$,
and there exists $X\subseteq V(G)$
with $|X|\ge \delta|G|$, such that one of $\overline{G}[X],G[X]$ has maximum degree at most $\vare\delta|G|$. (We are not free to replace $G$ by its complement, since the class of $\mathcal{H}$-free graphs is not closed under taking complements.)

Suppose that $\overline{G}[X]$ has maximum degree at most $\vare\delta|G|$. By \ref{pure1weak} applied to $\overline{G}$, 
there exist disjoint $A,B\subseteq X$, with $A$ complete to $B$, and with 
$|A|,|B|\ge \vare\delta|G|$. By \ref{parties}, $ \vare\delta|G|<|G|2^{-1/\tau}$, and so $ \vare\delta<2^{-1/\tau}$, contrary 
to the choice of $\tau$.

Thus $G[X]$ has maximum degree at most $\vare\delta|G|$.
Exactly as in the proof of \ref{treeext}, we obtain the blockade $\mathcal{B}$, and prove 
there is a $\mathcal{B}$-rainbow copy of $H$ or of $\overline{H}$.
The second is impossible since $G$ is $\overline{H}$-free; and so $G$ contains the star-expansion of $H$.
This proves \ref{antiforest}.~\bbox
 
We see that \ref{antiforestandhole} follows from \ref{antiforest}, by applying \ref{antiforest} to a forest $H'$ containing $H$ and containing 
a path of length at least 
$|E(C)|-4$
(because then the star-expansion of $H'$ contains $C$.)
The following is a theorem of  Bonamy, Bousquet and Thomass\'e~\cite{bonamy}:
\begin{thm}\label{bonamy}
For every integer $\ell>0$, there exists $\vare>0$ such that if $G$ has maximum degree less than $\vare|G|$, and $G$
has no anticomplete pair $(A,B)$ with $|A|,|B|\ge \vare|G|$, then $G$ has a hole of length at least $\ell$.
\end{thm}
The proof of \ref{antiforest} can be modified to show the following, by using \ref{bonamy} in place of \ref{pure1weak} (we omit the details):
\begin{thm}\label{antiholesandtree}
Let $H$ be the star-expansion of a forest; then for every integer $\ell\ge 3$,
$$\{H,\overline{C_{\ell}}, \overline{C_{\ell+1}},\overline{C_{\ell+2}},\ldots\}$$ 
has the Erd\H os-Hajnal property.
\end{thm}

This implies \ref{antiholesandhole}, by letting $H$ be the star-expansion of a path of length $|E(C)|-4$. (We may assume that
$C$ has length at least five, because it is known that $C_3,C_4$ both have the Erd\H os-Hajnal property.)

\section{$C_5$ with a hat}

There is still one result mentioned in the introduction that is not contained in any of the results we proved so far, namely
\ref{C5hatthm}, and now we will prove that.
\begin{thm}\label{C5hatthm2}
$\mathcal{H}=\{\widehat{C_5},\overline{\widehat{C_5}}\}$ has the Erd\H os-Hajnal property.
\end{thm}
\Proof
We proceed as usual: as in all these proofs, we choose a suitable $\vare\le 1/20$, choose $\delta$ satisfying \ref{rodl},
and then choose $\tau>0$ satisfying \ref{usesparse}, and we can also make $\tau$ less than
any positive function of the other parameters we choose. Let us see what we need.

We may assume (for a contradiction) that there is a $\tau$-critical $\mathcal{H}$-free graph $G$;
and there exists $X\subseteq V(G)$
with $|X|\ge \delta|G|$, such that $G[X]$ has maximum degree at most $\vare\delta|G|$. (We can pass to the complement if necessary.)
Let $\gamma=\delta/(400\vare)$.
By \ref{usesparse}, 
there is
a $(t,\gamma|G|/t^2)$-comb $((a_i,B_i):1\le i\le t)$ of $G[X]$ such that $t\ge 1/(400\vare)$ and $\{a_1\ll a_t\}$ is stable, and
there is a vertex $v\in X$ adjacent to $a_1\ll a_t$ and with
no neighbours in $B_1\cupcup B_t$.
\\
\\
(1) {\em For $1\le i\le t$, there is a component $G[D_i]$ of $G[B_i]$ with $|D_i|\ge \gamma|G|/t^3$.}
\\
\\
Suppose not, say for $i = 1$. 
Choose $s$ maximum such that there are $s$ subgraphs $F_1\ll F_s$
of $G[B_1]$, pairwise disjoint, each a union of components of $G[B_1]$, and each with at least $\gamma|G|/t^3$ vertices. We may assume 
that each $F_j$ is minimal,
and so has at most $2\gamma|G|/t^3$ vertices, since each component of $G[B_1]$ has at most $\gamma|G|/t^3$ vertices. Thus
$F_1\cupcup F_s$ has at most $2s\gamma|G|/t^3$ vertices, and so there are at least $\gamma|G|/t^2-2s\gamma|G|/t^3$ vertices of $B_1$
not in any of $F_1\ll F_s$. From the maximality of $s$, $\gamma|G|/t^2-2s\gamma|G|/t^3<\gamma|G|/t^3$, and so
$t-2s<1$. Hence $s\ge t/2$. But this contradicts \ref{parties}, since
we will arrange that $\gamma|G|/t^3\ge |G|(t/2)^{-1/\tau}$. To ensure this last, arrange at the start of the proof that 
$t^{1-3\tau}\ge 4$, by choosing $1/(400\vare)\ge 16$ and $\tau\le  1/6$, and arrange that $\gamma^\tau\ge 1/2$, by
choosing $\tau$ very small. This proves (1).
\\
\\
(2) {\em $\mathcal{D}=(D_1\ll D_t)$ is a pure blockade.}
\\
\\
Suppose not; then there exist distinct $i,j\in \{1\ll t\}$, such that some vertex $u\in D_j$ has both a neighbour and a 
non-neighbour in $D_i$.
Since $G[D_i]$ is connected, there is an edge $xy$ of $G[D_i]$ such that $u$ is adjacent 
to $x$ and not to $y$; and then the subgraph induced on $\{v,a_i, a_j, x,y,u\}$ is isomorphic to $\widehat{C_5}$, a contradiction. 
This proves (2).
\\
\\
(3) {\em There is no $\mathcal{D}$-rainbow triangle.}
\\
\\
Suppose there is, 
and so $G$ contains the star-expansion of $K_3$; but the star-expansion of $K_3$
contains $\widehat{C_5}$, a contradiction. This proves (3).

\bigskip

Let $P$ be the pattern of the pure blockade $\mathcal{D}$. Since $P$ is triangle-free by (3),  and $|P|=t$, it follows that there is
a stable set $I$  of $P$ with cardinality at least $t^{1/2}/2$. Hence the sets $D_i\;(i\in I)$ are pairwise anticomplete, but
we will arrange that $\gamma|G|/t^3\ge |G|(t^{1/2}/2)^{-1/\tau}$, a contradiction to \ref{parties}. To ensure this last, we arrange at the start of the proof that
$1/(400\vare)\ge 256$ and $\tau\le 1/12$, implying that $t\ge 256$ and $t^{1/2-3\tau}\ge 4$; and arrange that $\gamma^\tau\ge 1/2$ by choosing 
$\tau$ sufficiently small. This proves \ref{C5hatthm}.~\bbox

\end{document}